\documentclass[12pt,a4paper]{article}
\usepackage{amsmath}
\usepackage{amssymb}
\usepackage{t1enc}
\usepackage[latin1]{inputenc}
\usepackage[english]{babel}
\usepackage{amsfonts}
\usepackage{latexsym}
\usepackage{bm}
\newtheorem{theorem}{Theorem}[section]

\newtheorem{lemma}[theorem]{Lemma}

\begin{document}
\title{The maximum product of sizes of cross-$t$-intersecting uniform families}

\author{Peter Borg\\[5mm]
Department of Mathematics, University of Malta, Malta\\
%\texttt{p.borg.02@cantab.net}
peter.borg@um.edu.mt}
\date{}
\maketitle
%\date{13th June 2010}

\begin{abstract}
We say that a set $A$ \emph{$t$-intersects} a set $B$ if $A$ and $B$ have at least $t$ common elements. Two families $\mathcal{A}$ and $\mathcal{B}$ are said to be \emph{cross-$t$-intersecting} if each set in $\mathcal{A}$ $t$-intersects each set in $\mathcal{B}$. For any positive integers $n$ and $r$, let ${[n] \choose r}$ denote the family of all $r$-element subsets of $\{1,2,\dots, n\}$. We show that for any integers $r$, $s$ and $t$ with $1 \leq t \leq r \leq s$, there exists an integer $n_0(r,s,t)$ such that for any integer $n \geq n_0(r,s,t)$, if $\mathcal{A} \subset {[n] \choose r}$ and $\mathcal{B} \subset {[n] \choose s}$ such that $\mathcal{A}$ and $\mathcal{B}$ are cross-$t$-intersecting, then $|\mathcal{A}||\mathcal{B}| \leq {n-t \choose r-t}{n-t \choose s-t}$, and equality holds if and only if for some $T \in {[n] \choose t}$, $\mathcal{A} = \{A \in {[n] \choose r} \colon T \subset A\}$ and $\mathcal{B} = \{B \in {[n] \choose s} \colon T \subset B\}$.  This verifies a conjecture of Hirschorn.
\end{abstract}

\section{Introduction}\label{Intro}
%Give the parameter $\alpha_{\mathcal{F}}$ some suitable name.
%Maybe "cross-intersection parameter"?

Unless otherwise stated, throughout this paper we shall use small
letters such as $x$ to denote positive integers or elements of a set, capital letters such as $X$ to denote sets, and calligraphic letters such as $\mathcal{F}$ to denote \emph{families} (that is, sets whose elements are sets themselves). Arbitrary sets and families are taken to be finite and may be the \emph{empty set} $\emptyset$. An \emph{$r$-set} is a set of size $r$, that is, a set having
exactly $r$ elements. For any $n \geq 1$, $[n]$ denotes the set $\{1, \dots, n\}$ of the first $n$ positive integers. For any set $X$, ${X \choose r}$ denotes the set $\{A \subset X \colon |A| = r\}$ of all $r$-subsets of $X$. For any family $\mathcal{F}$, we denote the family $\{F \in \mathcal{F} \colon |F| = r\}$ by $\mathcal{F}^{(r)}$, and for any $t$-set $T$, we denote the family $\{F \in \mathcal{F} \colon T \subseteq F\}$ by $\mathcal{F}[T]$, and we call $\mathcal{F}[T]$ a \emph{$t$-star of $\mathcal{F}$} if $\mathcal{F}[T] \neq \emptyset$. 

Given an integer $t \geq 1$, we say that a set $A$
\emph{$t$-intersects} a set $B$ if $A$ and $B$ have at least $t$
common elements. A family $\mathcal{A}$ is said to be
\emph{$t$-intersecting} if each set in $\mathcal{A}$
$t$-intersects all the other sets in $\mathcal{A}$ (i.e.~$|A \cap B|
\geq t$ for every $A, B \in \mathcal{A}$ with $A \neq B$). A
$1$-intersecting family is also simply called an
\emph{intersecting family}. Note that $t$-stars are $t$-intersecting families.

Families $\mathcal{A}_1, \dots, \mathcal{A}_k$ are said to be \emph{cross-$t$-intersecting} if
for every $i$ and $j$ in $[k]$ with $i \neq j$, each set in $\mathcal{A}_i$ $t$-intersects each set in $\mathcal{A}_j$ (i.e.~$|A \cap B| \geq t$ for every $A \in \mathcal{A}_i$ and every $B \in \mathcal{A}_j$). Cross-$1$-intersecting families are also simply called
\emph{cross-intersecting families}.

The study of intersecting families took off with \cite{EKR}, which features the classical result that says that if $r \leq n/2$, then the size of a largest intersecting subfamily of ${[n] \choose
r}$ is ${n-1 \choose r-1}$, which is the size of every $1$-star of ${[n] \choose r}$. This result is known as the Erd\H os-Ko-Rado (EKR) Theorem. There are various proofs of this theorem (see \cite{HM,K,D}), two of which are particularly short and beautiful: Katona's \cite{K}, introducing the elegant cycle method, and Daykin's \cite{D}, using the powerful Kruskal-Katona Theorem \cite{Ka,Kr}. Also in \cite{EKR}, Erd\H os, Ko and Rado proved that for $t \leq r$, there exists an integer $n_0(r,t)$ such that for any $n \geq n_0(r,t)$, the size of a largest $t$-intersecting subfamily of ${[n] \choose r}$ is ${n-t \choose r-t}$, which is the size of every $t$-star of ${[n] \choose r}$. Frankl \cite{F_t1} showed that for $t \geq 15$, the smallest such $n_0(r,t)$ is $(r-t+1)(t+1)$. Subsequently, Wilson \cite{W} proved this for all $t \geq 1$. Frankl \cite{F_t1} conjectured that the size of a largest $t$-intersecting subfamily of ${[n] \choose r}$ is $\max\{|\{A \in {[n] \choose r} \colon |A \cap [t+2i]| \geq t+i\}| \colon i \in \{0\} \cup [r-t]\}$. A remarkable proof of this conjecture together with a complete characterisation of the extremal structures was obtained by Ahlswede and Khachatrian \cite{AK1}. 
%The classical Erd\H os-Ko-Rado (EKR) Theorem \cite{EKR} says that if %$n$ is sufficiently larger than $r$, then the size of any $t$-%intersecting subfamily of ${[n] \choose r}$ is at most ${n-t %\choose r-t}$, which is the number of sets
%in the $t$-intersecting subfamily of ${[n] \choose r}$
%consisting of those sets having $[t]$ as a subset. 
The $t$-intersection problem for $2^{[n]}$ was completely solved in \cite{Kat}. These are prominent results in extremal set theory. The EKR Theorem inspired a wealth of results of this kind, that is, results that establish how large a system of sets can be under certain
intersection conditions; see \cite{DF,F,F2,Borg7}.

For $t$-intersecting subfamilies of a given family $\mathcal{F}$,
the natural question to ask is how large they can be. For
cross-$t$-intersecting families, two natural parameters arise: the
sum and the product of sizes of the cross-$t$-intersecting
families (note that the product of sizes of $k$ families
$\mathcal{A}_1, \dots, \mathcal{A}_k$ is the number of $k$-tuples
$(A_1, \dots, A_k)$ such that $A_i \in \mathcal{A}_i$ for each $i
\in [k]$). It is therefore natural to consider the problem of
maximising the sum or the product of sizes of $k$
cross-$t$-intersecting subfamilies (not necessarily distinct or
non-empty) of a given family $\mathcal{F}$. The paper \cite{Borg8} analyses this problem in general and reduces it to the $t$-intersection problem for $k$ sufficiently large. %, but also outlines many of the known solutions for important families $\mathcal{F}$. 
In this paper we are concerned with the family ${[n] \choose r}$. We point out that the maximum product problem for $2^{[n]}$ and $k=2$ is solved in \cite{MT2}, and the maximum sum problem for $2^{[n]}$ and any $k$ is solved in \cite{Borg8} via the results in \cite{Kat,Kl,WZ}.
%results that solve such a problem for some particular important family %$\mathcal{F}$.

Wang and Zhang \cite{WZ} solved the maximum sum problem for ${[n] \choose r}$ using an elegant combination of the method used in \cite{Borg4,Borg3,Borg2,Borg5} and an important lemma that is found in \cite{AC,CK} and referred to as the `no-homomorphism lemma'. The solution for the case $t = 1$ had been obtained by Hilton \cite{H} and is the first result of this kind.

The maximum product problem for ${[n] \choose r}$ was first addressed by Pyber \cite{Pyber}, who proved that for any $r$, $s$ and $n$ such that either $r = s \leq n/2$ or $r < s$ and $n \geq 2s + r -2$, if $\mathcal{A} \subset {[n] \choose r}$ and $\mathcal{B} \subset {[n] \choose s}$ such that $\mathcal{A}$ and $\mathcal{B}$ are cross-intersecting, then $|\mathcal{A}||\mathcal{B}| \leq {n-1 \choose r-1}{n-1 \choose s-1}$. Subsequently, Matsumoto and Tokushige \cite{MT} proved this for any $r \leq s \leq n/2$, and they also determined the optimal structures. This brings us to the result of this paper, which solves the cross-$t$-intersection problem for $n$ sufficiently large and hence verifies \cite[Conjecture~3]{Hirschorn}.
%
%Hirschorn conjectured that if $n$ is sufficiently large, then $|%\mathcal{A}||\mathcal{B}| \leq {n-t \choose r-t}{n-t \choose s-t}$ for %every $\mathcal{A} \subset {[n] \choose r}$ and $\mathcal{B} \subset %{[n] \choose s}$ such that $\mathcal{A}$ and $\mathcal{B}$ are %cross-$t$-intersecting. We prove this conjecture and we also determine %the optimal structures. 
For $t \leq r \leq s$, let 
\[n_0(r,s,t) = \max\left\{r(s-t) {r+s-t \choose t}, \, (r-t){r \choose t} {r+s-t \choose t+1}\right\} + t+1.\]
\begin{theorem} \label{main} Let $t \leq r \leq s$ and $n \geq n_0(r,s,t)$. If  $\mathcal{A} \subset {[n] \choose r}$ and $\mathcal{B} \subset {[n] \choose s}$ such that $\mathcal{A}$ and $\mathcal{B}$ are cross-$t$-intersecting, then 
\[|\mathcal{A}||\mathcal{B}| \leq {n-t \choose r-t}{n-t \choose s-t},\]
and equality holds if and only if for some $T \in {[n] \choose t}$, $\mathcal{A} = \{A \in {[n] \choose r} \colon T \subset A\}$ and $\mathcal{B} = \{B \in {[n] \choose s} \colon T \subset B\}$.
\end{theorem}
The special case $r = s$ is treated in \cite{Tok1,Tok2,FLST}, which establish significantly better values of $n_0(r,r,t)$ that are close to optimal.

Theorem~\ref{main} generalises to one for $k \geq 2$ cross-$t$-intersecting families.
\begin{theorem} \label{maincor} Let $k \geq 2$, $t \leq r_1 \leq \dots \leq r_k$ and $n \geq n_0(r_{k-1},r_k,t)$. If $\mathcal{A}_1 \subset {[n] \choose r_1}, \dots, \mathcal{A}_k \subset {[n] \choose r_k}$, and $\mathcal{A}_1, \dots, \mathcal{A}_k$ are cross-$t$-intersecting, then 
\[\prod_{i=1}^k |\mathcal{A}_i| \leq \prod_{i=1}^k {n-t \choose r_i - t},\]
and equality holds if and only if for some $T \in {[n] \choose t}$, $\mathcal{A}_i = \{A \in {[n] \choose r_i} \colon T \subset A\}$ for each $i \in [k]$.
\end{theorem}
\textbf{Proof.} For each $i \in [k]$, let $a_i = |\mathcal{A}_i|$ and $b_i = {n-t \choose r_i-t}$. Note that $n_0(r_i,r_j,t) \leq n_0(r_{k-1},r_k,t)$ for any $i, j \in [k]$ with $i \neq j$. Thus, by Theorem~\ref{main}, for any $i, j \in [k]$ with $i \neq j$, we have $a_ia_j \leq b_ib_j$, and equality holds if and only if for some $T_{i,j} \in {[n] \choose t}$, $\mathcal{A}_i = \{A \in {[n] \choose r_i} \colon T_{i,j} \subset A\}$ and $\mathcal{A}_j = \{A \in {[n] \choose r_j} \colon T_{i,j} \subset A\}$. Let mod$^*$ represent the usual \emph{modulo operation} with the exception that for any two integers $x$ and $y > 0$, $(xy) \, {\rm mod}^* \, y$ is $y$ instead of $0$. We have
\begin{align} \left( \prod_{i=1}^k a_i \right)^2 &=
(a_1a_2)(a_{3 \, {\rm mod}^* \, k}a_{4 \, {\rm mod}^* \,
k})\cdots(a_{(2k-1) \, {\rm mod}^* \, k} a_{(2k) \, {\rm mod}^*
\, k}) \nonumber \\
&\leq (b_1b_2)(b_{3 \, {\rm mod}^* \, k}b_{4 \, {\rm mod}^* \,
k})\cdots(b_{(2k-1) \, {\rm mod}^* \, k} b_{(2k) \, {\rm mod}^*
\, k}) = \left( \prod_{i=1}^k b_i \right)^2. \nonumber
\end{align}
So $\prod_{i=1}^k a_i \leq \prod_{i=1}^k b_i$. Suppose equality holds. Then for any $i, j \in [k]$ with $i \neq j$, $\mathcal{A}_i = \{A \in {[n] \choose r_i} \colon T_{i,j} \subset A\}$ and $\mathcal{A}_j = \{A \in {[n] \choose r_j} \colon T_{i,j} \subset A\}$ for some $T_{i,j} \in {[n] \choose t}$. So we have $\{A \in {[n] \choose r_1} \colon T_{1,2} \subset A\} = \mathcal{A}_1 = \{A \in {[n] \choose r_1} \colon T_{1,j} \subset A\}$ for each $j \in [k] \backslash \{1\}$. So $T_{1,2} = T_{1,j}$ for each $j \in [k] \backslash \{1\}$. So $\mathcal{A}_j = \{A \in {[n] \choose r_j} \colon T_{1,2} \subset A\}$ for each $j \in [k] \backslash \{1\}$. Hence the result.~\hfill{$\Box$}

\section{The compression operation}
\label{Compsection}
For any $i, j \in [n]$, let $\delta_{i,j}
\colon 2^{[n]} \rightarrow 2^{[n]}$ be defined by
\[ \delta_{i,j}(A) = \left\{ \begin{array}{ll}
(A \backslash \{j\}) \cup \{i\} & \mbox{if $j \in A$ and $i \notin
A$};\\
A & \mbox{otherwise,}
\end{array} \right. \]
and let $\Delta_{i,j} \colon 2^{2^{[n]}} \rightarrow 2^{2^{[n]}}$ be the
\emph{compression operation} (see \cite{EKR}) defined by
\[\Delta_{i,j}(\mathcal{A}) = \{\delta_{i,j}(A) \colon A \in
\mathcal{A}, \delta_{i,j}(A) \notin \mathcal{A}\} \cup \{A \in
\mathcal{A} \colon \delta_{i,j}(A) \in \mathcal{A}\}.\]
Note that $|\Delta_{i,j}(\mathcal{A})| = |\mathcal{A}|$. \cite{F}
provides a survey on the properties and uses of compression (also
called \emph{shifting}) operations in extremal set theory. %We will need the following basic result, which we prove for completeness.

If $i < j$, then we call $\Delta_{i,j}$ a \emph{left-compression}. A family $\mathcal{F} \subseteq 2^{[n]}$ is said to be \emph{compressed} if $\Delta_{i,j}(\mathcal{F}) = \mathcal{F}$ for every $i,j \in [n]$ with $i < j$. In other words, $\mathcal{F}$ is compressed if it is invariant under left-compressions.

The following lemma captures some well-known fundamental properties of compressions, and we will prove it for completeness.

\begin{lemma}\label{compcross} Let $\mathcal{A}$ and $\mathcal{B}$ be cross-$t$-intersecting subfamilies of $2^{[n]}$.\\
(i) For any $i, j \in [n]$, $\Delta_{i,j}(\mathcal{A})$ and $\Delta_{i,j}(\mathcal{B})$ are cross-$t$-intersecting subfamilies of $2^{[n]}$.\\
(ii) If $t \leq r \leq s \leq n$, $\mathcal{A} \subset {[n] \choose r}$, $\mathcal{B} \subset {[n] \choose s}$, and $\mathcal{A}$ and $\mathcal{B}$ are compressed, then $|A \cap B \cap [r+s-t]| \geq t$ for any $A \in \mathcal{A}$ and any $B \in \mathcal{B}$.
\end{lemma}
\textbf{Proof.} %The result is trivial if $\Delta_{i,j}(\mathcal{A}) = %\emptyset$ or $\Delta_{i,j}(\mathcal{B}) = \emptyset$.
%%Consider $\Delta_{i,j}(\mathcal{A}) \neq \emptyset$ and $\Delta_{i,j}
%%(\mathcal{B}) \neq \emptyset$.
Let $i, j \in [n]$. Suppose $A \in \Delta_{i,j}(\mathcal{A})$ and $B \in \Delta_{i,j}(\mathcal{B})$. If $A \in \mathcal{A}$ and $B \in \mathcal{B}$, then $|A \cap B| \geq t$ since $\mathcal{A}$ and $\mathcal{B}$ are cross-$t$-intersecting. Suppose $A \notin \mathcal{A}$ or $B \notin \mathcal{B}$; we may assume that $A \notin \mathcal{A}$. Then $A = \delta_{i,j}(A') \neq A'$ for some $A' \in \mathcal{A}$.  So $i \notin A'$, $j \in A'$, $i \in A$ and $j \notin A$. Suppose $|A \cap B| \leq t-1$. Then $i \notin B$ and hence $B \in \mathcal{B}$. So $B \in \mathcal{B} \cap \Delta_{i,j}(\mathcal{B})$ and hence $B, \delta_{i,j}(B) \in \mathcal{B}$. So $|A' \cap B| \geq t$ and $|A' \cap \delta_{i,j}(B)| \geq t$. From $|A \cap B| \leq t-1$ and $|A' \cap B| \geq t$ we get $A' \cap B = (A \cap B) \cup \{j\}$ and hence $A' \cap \delta_{i,j}(B) = A \cap B$, but this contradicts $|A \cap B| \leq t-1$ and $|A' \cap \delta_{i,j}(B)| \geq t$. Hence (i).

Suppose $t \leq r \leq s \leq n$, $\mathcal{A} \subset {[n] \choose r}$ and $\mathcal{B} \subset {[n] \choose s}$. Let $A \in \mathcal{A}$ and $B \in \mathcal{B}$. So $|A \cap B| \geq t$. Let $X := (A \cap B) \cap [r+s-t]$, $Y = (A \cap B) \backslash [r+s-t]$ and $Z = [r+s-t] \backslash (A \cup B)$. If $Y = \emptyset$, then $X = A \cap B$ and hence $|X| \geq t$. Now consider $Y \neq \emptyset$. Let $p =
|Y|$. Since
\begin{align} |Z| &= r + s - t - |(A \cup B) \cap [r+s-t]| \geq r + s - t - |X| - |A \backslash B| - |B \backslash A| \nonumber \\
&= r + s - t - |X| - |A \backslash (X \cup Y)| - |B
\backslash (X \cup Y)| \nonumber \\
&= r+s-t - |X| - (r - |X| - |Y|) - (s - |X| - |Y|) = 2|Y| + |X| - t \nonumber \\
&= |Y| + |Y \cup X| - t = p + |A \cap B| - t \geq p, \nonumber
\end{align}
${Z \choose p} \neq \emptyset$. Let $W \in {Z \choose p}$. Let $C
:= (B \backslash Y) \cup W$. Let $y_1, \dots, y_p$ be the elements
of $Y$, and let $w_1, \dots, w_p$ be those of $W$. So $C =
\delta_{w_1,y_1} \circ \dots \circ \delta_{w_p,y_p}(B)$. Note that
$\delta_{w_1,y_1}, \dots, \delta_{w_p,y_p}$ are left-compressions
as $W \subseteq [r+s-t]$ and $Y \subseteq [n] \backslash [r+s-t]$.
Since $\mathcal{B}$ is left-compressed, $C \in \mathcal{B}$. So
$|A \cap C| \geq t$. Now clearly $|A \cap C| = |X|$. So $|X| \geq t$. Hence (ii).~\hfill{$\Box$}\\

%Suppose not. Then there exists $A \in \Delta_{i,j}(\mathcal{A})$ and $B %\in \Delta_{i,j}(\mathcal{B})$ such that $A \cap B = \emptyset$. Since %$\mathcal{A}$ and $\mathcal{B}$ are cross-intersecting, we cannot have %both $A \in \mathcal{A}$ and $B \in \mathcal{B}$. We may assume that $A %\notin \mathcal{A}$. So $A = \delta_{i,j}(C) \neq C$ for some $C \in %\mathcal{A}$. So $j \in C$ and $i \notin C$. $A' = A \backslash \{i\} %\cup \{j\}$

Suppose a subfamily $\mathcal{A}$ of $2^{[n]}$ is not compressed. Then $\mathcal{A}$ can be transformed to a compressed family through left-compressions as follows. Since $\mathcal{A}$ is not compressed, we can find a left-compression that changes $\mathcal{A}$, and we apply it to $\mathcal{A}$ to obtain a new subfamily of $2^{[n]}$. We keep on repeating this (always applying a left-compression to the last family obtained) until we obtain a subfamily of $2^{[n]}$ that is invariant under every left-compression (such a point is indeed reached, because if $\Delta_{i,j}(\mathcal{F}) \neq \mathcal{F} \subseteq 2^{[n]}$ and $i < j$, then $0 < \sum_{G \in \Delta_{i,j}(\mathcal{F})} \sum_{b \in G} b < \sum_{F \in \mathcal{F}} \sum_{a \in F} a$).

Now consider $\mathcal{A}, \mathcal{B} \subseteq 2^{[n]}$ such that $\mathcal{A}$ and $\mathcal{B}$ are cross-$t$-intersecting. Then, by Lemma~\ref{compcross}(i), we can obtain $\mathcal{A}^*, \mathcal{B}^* \subseteq 2^{[n]}$ such that $\mathcal{A}^*$ and $\mathcal{B}^*$ are compressed and cross-$t$-intersecting, $|\mathcal{A}^*| = |\mathcal{A}|$ and $|\mathcal{B}^*| = |\mathcal{B}|$. Indeed, similarly to the above procedure, if we can find a left-compression that changes at least one of $\mathcal{A}$ and $\mathcal{B}$, then we apply it to both $\mathcal{A}$ and $\mathcal{B}$, and we keep on repeating this (always performing this on the last two families obtained) until we obtain $\mathcal{A}^*, \mathcal{B}^* \subseteq 2^{[n]}$ such that both $\mathcal{A}^*$ and $\mathcal{B}^*$ are invariant under every left-compression.

\section{Proof of Theorem~\ref{main}} \label{Proof}

%\begin{lemma}\label{int comp} Let $t \leq r \leq s$ and $n \geq r+s-t$. %Let $\mathcal{A}$ and $\mathcal{B}$ be cross-$t$-intersecting left-%compressed subfamilies of $2^{[n]}$ such that $\emptyset \neq %\mathcal{A} \subset {[n] \choose r}$ and $\emptyset \neq \mathcal{B} %\subset {[n] \choose s}$. Then $|A \cap B \cap [r+s-t]| \geq t$ for any %$A \in \mathcal{A}$ and $B \in \mathcal{B}$.
%\end{lemma}
%

We will need the following lemma only when dealing with the characterisation of the extremal structures in the proof of Theorem~\ref{main}. 

\begin{lemma}\label{charlemma} Let $r$, $s$, $t$ and $n$ be as in Theorem~\ref{main}, and let $i,j \in [n]$. Let $\mathcal{H} = 2^{[n]}$. Let $\mathcal{A} \subset \mathcal{H}^{(r)}$ and $\mathcal{B} \subset \mathcal{H}^{(s)}$ such that $\mathcal{A}$ and $\mathcal{B}$ are cross-$t$-intersecting. Suppose $\Delta_{i,j}(\mathcal{A}) = \mathcal{H}^{(r)}[T]$ and $\Delta_{i,j}(\mathcal{B}) = \mathcal{H}^{(s)}[T]$ for some $T \in {[n] \choose t}$. Then $\mathcal{A} = \mathcal{H}^{(r)}[T']$ and $\mathcal{B} = \mathcal{H}^{(s)}[T']$ for some $T' \in {[n] \choose t}$. 
\end{lemma}
We prove this lemma using the following special case of \cite[Lemma~5.6]{Borg}. 
%%
%\begin{lemma}[{\cite[Lemma~5.6]{Borg}}] \label{L comp} Let $t \leq q$, %$P \subseteq [t,q]$ and $n \geq 2q-t+1$. Let $\mathcal{H} = 2^{[n]}$. %Let $\mathcal{G}$ be a
%$t$-intersecting subfamily of $\bigcup_{p \in P}
%\mathcal{H}^{(p)}$. Suppose $\Delta_{i,j}(\mathcal{G})$ is a largest %$t$-star of $\bigcup_{p \in P} \mathcal{H}^{(p)}$ for some $i, j \in
%[n]$ with $i \neq j$. Then $\mathcal{G}$ is a largest $t$-star of %$\bigcup_{p \in P} \mathcal{H}^{(p)}$.
%\end{lemma}
%%
%
\begin{lemma} \label{L comp} Let $t \leq p$ and $n \geq 2p-t+1$. Let $\mathcal{H} = 2^{[n]}$. Let $\mathcal{G}$ be a $t$-intersecting subfamily of $\mathcal{H}^{(p)}$. Suppose $\Delta_{i,j}(\mathcal{G})$ is a largest $t$-star of $\mathcal{H}^{(p)}$ for some $i, j \in [n]$. Then $\mathcal{G}$ is a largest $t$-star of $\mathcal{H}^{(p)}$.
\end{lemma}
\textbf{Proof of Lemma~\ref{charlemma}.} Note that $T \backslash \{i\} \subset E$ for all $E \in \mathcal{A} \cup \mathcal{B}$. 

Suppose $\mathcal{A}$ is not $t$-intersecting. Then there exist $A_1, A_2 \in \mathcal{A}$ such that $|A_1 \cap A_2| \leq t-1$. So $T \nsubseteq A_l$ for some $l \in [2]$; we may (and will) assume that $l=1$. Thus, since $\Delta_{i,j}(\mathcal{A}) = \mathcal{H}^{(r)}[T]$, we must have $A_1 \neq \delta_{i,j}(A_1) \in \Delta_{i,j}(\mathcal{A})$, $\delta_{i,j}(A_1) \notin \mathcal{A}$ (because otherwise $A_1 \in \Delta_{i,j}(\mathcal{A})$), $i \in T$, $j \notin T$, $j \in A_1$ and $A_1 \cap T = T \backslash \{i\}$. Since $T \backslash \{i\} \subset A_1 \cap A_2$ and $|A_1 \cap A_2| \leq t-1$, we have $A_1 \cap A_2 = T \backslash \{i\}$. So $j \notin A_2$ and  
%$A_2 \neq \delta_{i,j}(A_2)$; then $j \in A_2$, $T \backslash \{i\} %\in A_2$, and hence $|A_1 \cap A_2| \geq t$, a contradiction. 
hence $A_2 = \delta_{i,j}(A_2)$. Since $\delta_{i,j}(A_2) \in \Delta_{i,j}(\mathcal{A})$, $T \subseteq A_2$. 
%Since $|A_1 \cap A_2| \leq t-1$, we have $A_1 \cap A_2 = T \backslash %\{i\}$. 
Let $X = [n] \backslash (A_1 \cup A_2)$. Since $|X| = n - |A_1 \cup A_2| = n - (|A_1| + |A_2| - |A_1 \cap A_2|) = n - 2r + t \geq n_0(r,s,t) - 2r + t > s-t$, ${X \choose s-t} \neq \emptyset$. Let $C \in {X \choose s-t}$ and $D = C \cup T$. So $D \in \mathcal{H}^{(s)}[T]$ and $D \cap A_1 = T \backslash \{i\}$, meaning that $D \in \Delta_{i,j}(\mathcal{B})$ and $|D \cap A_1| = t-1$. Thus, since $\mathcal{A}$ and $\mathcal{B}$ are cross-$t$-intersecting, $D \notin \mathcal{B}$ and $(D \backslash \{i\}) \cup \{j\} \in \mathcal{B}$, but this is a contradiction since $|((D \backslash \{i\}) \cup \{j\}) \cap A_2| = |T \backslash \{i\}| = t-1$. 

Therefore, $\mathcal{A}$ is $t$-intersecting. Similarly, $\mathcal{B}$ is $t$-intersecting. %Since $\mathcal{A}$ and $\mathcal{B}$ are %cross-$t$-intersecting, it follows that $\mathcal{A} \cup \mathcal{B}$ %is $t$-intersecting. 
Now $\mathcal{H}^{(r)}[T]$ and $\mathcal{H}^{(s)}[T]$ are largest $t$-stars of $\mathcal{H}^{(r)}$ and $\mathcal{H}^{(s)}$, respectively. So $\Delta_{i,j}(\mathcal{A})$ and and $\Delta_{i,j}(\mathcal{B})$ are largest $t$-stars of $\mathcal{H}^{(r)}$ and $\mathcal{H}^{(s)}$, respectively. By Lemma~\ref{L comp}, for some $T', T^* \in {[n] \choose t}$, $\mathcal{A} = \mathcal{H}^{(r)}[T']$ and $\mathcal{B} = \mathcal{H}^{(s)}[T^*]$. Suppose $T' \neq T^*$. Then clearly we can find $A' \in \mathcal{H}^{(r)}[T']$ and $B' \in \mathcal{H}^{(s)}[T']$ such that $|A' \cap B'| \leq t-1$; however, this is a contradiction since $\mathcal{A}$ and $\mathcal{B}$ are cross-$t$-intersecting. So $T^* = T'$.~\hfill{$\Box$}\\
%So $\mathcal{H}^{(r)}[T] \cup \mathcal{H}^{(s)}[T]$ is a largest $t$-%star of $\mathcal{H}^{(r)} \cup \mathcal{H}^{(s)}$. 
%
%Suppose $r < s$. Then we have $\Delta_{i,j}(\mathcal{A} \cup %\mathcal{B}) = \Delta_{i,j}(\mathcal{A}) \cup \Delta_{i,j}%(\mathcal{B}) = \mathcal{H}^{(r)}[T] \cup %\mathcal{H}^{(s)}[T]$, and hence $\Delta_{i,j}(\mathcal{A} \cup %\mathcal{B})$ is a largest $t$-star of $\mathcal{H}^{(r)} \cup %\mathcal{H}^{(s)}$. By Lemma~\ref{L comp} (with $q = s$, $P = \{r,s\}$ %and $\mathcal{G} = \mathcal{A} \cup \mathcal{B}$), $\mathcal{A} \cup %\mathcal{B}$ is a largest $t$-star of $\mathcal{H}^{(r)} \cup %\mathcal{H}^{(s)}$. So $\mathcal{A} \cup \mathcal{B} = %(\mathcal{H}^{(r)} \cup \mathcal{H}^{(s)})[T'] = \mathcal{H}^{(r)}[T'] %\cup \mathcal{H}^{(s)}[T']$ for some $T' \in {[n] \choose t}$. So %$\mathcal{A} = \mathcal{H}^{(r)}[T']$ and $\mathcal{B} = %\mathcal{H}^{(s)}[T']$.
%
%Now suppose $r = s$. ~\hfill{$\Box$}\\ 
%Note that Lemma~\ref{charlemma} can be improved for a significant %range of values of $n$ that are smaller than $n_0(r,s,t)$; needing $n %\geq n_0(r,s,t)$ in the proof of Theorem~\ref{main} allowed us to %relax our argument for this lemma. \\
\\
\textbf{Proof of Theorem~\ref{main}.} Let $n_0 = n_0(r,s,t)$. Let $\mathcal{H} = 2^{[n]}$. So ${[n] \choose r} = \mathcal{H}^{(r)}$ and ${[n] \choose s} = \mathcal{H}^{(s)}$. If $\mathcal{A} = \emptyset$ or $\mathcal{B}= \emptyset$, then $|\mathcal{A}||\mathcal{B}| = 0$. So we assume that $\mathcal{A} \neq \emptyset$ and $\mathcal{B} \neq \emptyset$. 

%If $r = t$, then we trivially have that for some $T \in {[n] \choose %t}$, $\mathcal{A} = \{T\}$ and $\mathcal{B} \subseteq \mathcal{H}^{(s)}%[T]$. So we assume that $r \geq t+1$. 

As explained in Section~\ref{Compsection}, we apply left-compressions to $\mathcal{A}$ and $\mathcal{B}$ simultaneously until we obtain two compressed families $\mathcal{A}^*$ and $\mathcal{B}^*$, respectively, and we know that $\mathcal{A}^*$ and $\mathcal{B}^*$ are cross-$t$-intersecting, $\mathcal{A}^* \subset \mathcal{H}^{(r)}$, $\mathcal{B}^* \subset \mathcal{H}^{(s)}$, $|\mathcal{A}^*| = |\mathcal{A}|$ and $|\mathcal{B}^*| = |\mathcal{B}|$. In view of Lemma~\ref{charlemma}, we may therefore assume that $\mathcal{A}$ and $\mathcal{B}$ are compressed.

By Lemma~\ref{compcross}(ii), 
\begin{equation} \mbox{$|A \cap [r+s-t]| \geq t$ for each $A \in \mathcal{A}$.} \label{1}
\end{equation}

\textit{Case 1:} $|A^* \cap [r+s-t]| = t$ for some $A^* \in \mathcal{A}$. Then $A^* \cap [r+s-t] = T^*$ for some $T^* \in {[r+s-t] \choose t}$. By Lemma~\ref{compcross}(ii), $T^* \subset B$ for each $B \in \mathcal{B}$. So $\mathcal{B} = \mathcal{B}[T^*]$.

If $T^* \subset A$ for each $A \in \mathcal{A}$, then $|\mathcal{A}||\mathcal{B}| \leq |\mathcal{H}^{(r)}[T^*]||\mathcal{H}^{(s)}[T^*]| = {n-t \choose r-t}{n-t \choose s-t}$, and equality holds if and only if $\mathcal{A} = \mathcal{H}^{(r)}[T^*]$ and $\mathcal{B} = \mathcal{H}^{(s)}[T^*]$.

Suppose $T^* \nsubseteq A'$ for some $A' \in \mathcal{A}$. Then $|A' \cap T^*| \leq t-1$. Let $C = A' \cap T^*$ and $D = A' \backslash C$. For each $B \in \mathcal{B}$, we have $t \leq |B \cap A'| = |B \cap C| + |B \cap D| = |C| + |B \cap D| \leq t-1 + |B \cap D|$ and hence $|B \cap D| \geq 1$. So $\mathcal{B} = \{B \in \mathcal{B}[T^*] \colon |B \cap D| \geq 1\}$. Together with (\ref{1}), this gives us
\begin{align} |\mathcal{A}||\mathcal{B}| &= \left| \bigcup_{T \in {[r+s-t] \choose t}} \mathcal{A}[T] \right| \left| \bigcup_{X \in {D \choose 1}} \mathcal{B}[T^* \cup X]\right| \nonumber \\
&\leq \left(\sum_{T \in {[r+s-t] \choose t}} |\mathcal{A}[T]| \right) \left(\sum_{X \in {D \choose 1}} |\mathcal{B}[T^* \cup X]|\right) \nonumber \\
&\leq \left(\sum_{T \in {[r+s-t] \choose t}} |\mathcal{H}^{(r)}[T]| \right) \left(\sum_{X \in {D \choose 1}} |\mathcal{H}^{(s)}[T^* \cup X]|\right) \nonumber \\
&= \left(\sum_{T \in {[r+s-t] \choose t}} {n-t \choose r-t} \right) \left(\sum_{X \in {D \choose 1}} {n-t-1 \choose s-t-1}\right) \nonumber \\
&= {r+s-t \choose t} {n-t \choose r-t} {|D| \choose 1} {n-t-1 \choose s-t-1} \nonumber \\
&\leq r{r+s-t \choose t} {n-t \choose r-t} {n-t-1 \choose s-t-1} \nonumber \\
&= r{r+s-t \choose t} {n-t \choose r-t} \frac{s-t}{n-t}{n-t \choose s-t} \nonumber \\
&\leq \frac{s-t}{n_0-t}r{r+s-t \choose t} {n-t \choose r-t} {n-t \choose s-t} \nonumber \\
&< {n-t \choose r-t}{n-t \choose s-t}. \nonumber
%&= r{r+s-t \choose t} \prod_{i=0}^{r-t-1} \frac{n-t-i}{r-t-i} %\prod_{j=0}^{s-t-2} \frac{n-t-1-j}{s-t-1-j} \nonumber \\
%&\leq r{r+s-t \choose t} (n-r+1)^{r-t}(n-s+1)^{s-t-1} \nonumber \\
\end{align}

\textit{Case 2:} $|A \cap [r+s-t]| \geq t+1$ for all $A \in \mathcal{A}$. So $\mathcal{A} = \bigcup_{Z \in {[r+s-t] \choose t+1}}\mathcal{A}[Z] \subseteq \bigcup_{Z \in {[r+s-t] \choose t+1}}\mathcal{H}^{(r)}[Z]$. Let $A^* \in \mathcal{A}$. Since $|B \cap A^*| \geq t$ for all $B \in \mathcal{B}$, we have $\mathcal{B} = \bigcup_{T \in {A^* \choose t}}\mathcal{B}[T] \subseteq \bigcup_{T \in {A^* \choose t}}\mathcal{H}^{(s)}[T]$. Therefore,
\begin{align} %|\mathcal{A}||\mathcal{B}| &\leq \left(\sum_{Z \in %{[r+s-t] \choose t+1}} |\mathcal{H}^{(r)}[Z]| \right) \left(\sum_{T %\in {A^* \choose t}} |\mathcal{H}^{(s)}[T]|\right) \nonumber \\
%&= \left(\sum_{T \in {[r+s-t] \choose t}} {n-t \choose r-t} \right) %\left(\sum_{X \in {D \choose 1}} {n-t-1 \choose s-t-1}\right) %\nonumber \\
|\mathcal{A}||\mathcal{B}| &\leq {r+s-t \choose t+1} {n-t-1 \choose r-t-1} {r \choose t} {n-t \choose s-t} \nonumber \\
&= {r+s-t \choose t+1} \frac{r-t}{n-t}{n-t \choose r-t} {r \choose t} {n-t \choose s-t} \nonumber \\
&= \frac{r-t}{n_0-t}  {r \choose t} {r+s-t \choose t+1} {n-t \choose r-t} {n-t \choose s-t} \nonumber \\
&< {n-t \choose r-t}{n-t \choose s-t}. \nonumber
\end{align}
\hfill{$\Box$}

\end{document}